\documentclass[a4paper]{article}

\usepackage{a4wide}
\usepackage{amsmath}
\usepackage{amssymb} 
\usepackage{bm}
\usepackage{booktabs}
\usepackage{algorithm}
\usepackage{multicol}
\usepackage{multirow}
\usepackage[T1]{fontenc}
\usepackage{lmodern, textcomp}
\usepackage{amsfonts}
\usepackage{graphicx}
\usepackage{epstopdf}
\usepackage{algorithmic}
\usepackage{cleveref}
\usepackage{amsthm}
\usepackage{url}
\usepackage{color}
\usepackage{subcaption}

\theoremstyle{thmstyleone}%
\newtheorem{theorem}{Theorem}%
\newtheorem{lemma}{Lemma}%
\newtheorem{proposition}{Proposition}%

\theoremstyle{thmstyletwo}%

\newtheorem{remark}{Remark}

\theoremstyle{thmstylethree}

\title{Theoretical insights on the residual transformation from bi-conjugate gradient into bi-conjugate residual via\\a smoothing scheme}

\date{}

\author{Arisa Kawase\thanks{Graduate School of Integrative Science and Engineering, Tokyo City University, 1-28-1 Tamazutsumi, Setagaya-ku, Tokyo 158-8557, Japan.}
\and Kensuke Aihara\thanks{Department of Computer Science, Tokyo City University, 1-28-1 Tamazutsumi, Setagaya-ku, Tokyo, 158-8557, Japan ({\tt aiharak@tcu.ac.jp}).}}

\begin{document}

\maketitle

\begin{abstract}
Bi-conjugate gradient (Bi-CG) and bi-conjugate residual (Bi-CR) methods are underlying iterative solvers for linear systems with nonsymmetric matrices. 
Residual smoothing is a standard technique for obtaining smooth convergence behavior of residual norms; additionally, it represents the transformation between iterative methods.
For example, the residuals of the CR method can be obtained by applying a smoothing scheme to those of the CG method for symmetric linear systems. 
Based on this relationship, the transformation from Bi-CG residuals to Bi-CR residuals using a smoothing scheme was examined in our previous study [Kawase, A., Aihara, K.: Transformation from Bi-CG into Bi-CR Using a Residual Smoothing-like Scheme. AIP Conference Proceedings (2026)]; however, we only provided heuristic and experimental observations. 
In the present study, we provide a detailed discussion on the theoretical aspects of these transformations. 
Specifically, we prove that the resulting algorithm transformed from the Bi-CG method using the residual smoothing technique has the same bi-orthogonal properties as those of the original Bi-CR method. 
We also present a more concise transformation algorithm and its numerical example. 
These analyses complement our previous study and provide theoretical validity of the residual transformation between the Bi-CG and Bi-CR methods. 
\end{abstract}

\


\


\section{Introduction}\label{sec1}

In this study, we focus on Krylov subspace methods (e.g.,~\cite{Saad_2003,vanderVorst_2003}), which represent a framework of iterative methods for solving linear systems
\begin{align*}
A\bm{x} = \bm{b},\quad A\in \mathbb{R}^{n\times n},\quad \bm{b} \in \mathbb{R}^n, 
\end{align*}
where $A$ is a large and sparse nonsingular matrix. 

The conjugate gradient (CG) method~\cite{Hestenes_1952} is the most basic Krylov subspace method for symmetric positive definite (SPD) matrices, and the related conjugate residual (CR) method~\cite{Eisenstat_1983} is often used for symmetric (but not necessarily SPD) systems. 
The Bi-CG~\cite{Fletcher_1976} and Bi-CR~\cite{Sogabe_2005,Sogabe_2009} methods are known as underlying solvers that are extensions of the CG and CR methods, respectively, to nonsymmetric cases. 
In this study, we consider the relationships between the above methods; in particular, we investigate the transformation between the residuals generated by the Bi-CG and Bi-CR methods from theoretical viewpoints.

\subsection{Motivation and objective}

Residual smoothing techniques~\cite{Schonauer_1987,Weiss_1996} are useful for converting a sequence of residuals obtained by iterative methods into alternatives with smoother convergence behavior. 
Based on this property, a smoothing scheme represents a residual transformation between different iterative methods.
For example, applying minimal residual smoothing (MRS) to the CG method generates the same residuals as in the CR method~\cite{Walker_1995,Weiss_1996}. 
Analyzing this connection provides novel insights into iterative methods, enabling the development of new techniques for improving convergence. 
As a historically significant development, it has been demonstrated that the residuals of the Bi-CG method can be transformed into those generated by the quasi-minimal residual (QMR) method~\cite{Freund_1991} using a smoothing form, which leads to an alternative smoothing scheme (i.e., QMRS); see~\cite{Zhou_1994} for further details. 

The aforementioned studies motivated us to consider whether a similar residual transformation exists between the Bi-CG and Bi-CR methods. 
Thus, our previous study~\cite{Kawase_2025} derived a plausible transformation from Bi-CG residuals to Bi-CR residuals using a smoothing form. 
Although this derivation process is rational and the transformation has been validated through numerical experiments, its theoretical consistency has not been discussed. 
Therefore, the present study theoretically proves that the transformation algorithm \cite[Algorithm~4]{Kawase_2025} has bi-orthogonal properties for the residuals and direction vectors, which characterize the Bi-CR method. 
Specifically, the novelty of the present study lies in demonstrating the theoretical equivalence between \cite[Algorithm~4]{Kawase_2025} and the original Bi-CR method by proving the bi-orthogonal properties (described as Theorem~\ref{theorem} in subsection~\ref{sec_theorem}). 
This insight provides novel evidence that Bi-CR residuals can be generated via the Bi-CG method with a residual smoothing scheme. 
Moreover, we present a more concise algorithm for transforming the Bi-CG residuals into Bi-CR residuals. 
The algorithm, described as Algorithm~\ref{Bi-CR_3} in Sect.~\ref{sec5}, also highlights the novelty of this study over the previous study~\cite{Kawase_2025}. 
In addition, a numerical example is presented to support our theoretical analysis.

\subsection{Notation}

The Krylov subspace with $A$ and initial residual $\bm{r}_0 := \bm{b}-A\bm{x}_0$ is defined as
\begin{align*}
\mathcal K_k(A,\bm{r}_0) := \text{span}(\bm{r}_0, A\bm{r}_0,\dots, A^{k-1}\bm{r}_0), 
\end{align*}
where $\bm{x}_0 \in \mathbb{R}^n$ denotes the initial estimate. 
Below, the coefficient matrix $A$ is assumed to be nonsymmetric and nonsingular unless otherwise noted. 
Moreover, we only consider the case of $k < n$ using exact arithmetic, and a breakdown does not occur. 

The $H$-inner product is defined by $(\bm{x}, \bm{y})_H := \bm{x}^\top H\bm{y}$ for $\bm{x}, \bm{y} \in \mathbb{R}^{n}$ and an SPD matrix $H \in \mathbb{R}^{n\times n}$. 
The corresponding induced norm is denoted as $\|\cdot \|_H$. 
These are simplified to standard inner product $(\bm{x}, \bm{y}) = \bm{x}^\top \bm{y}$ and Euclidean norm $\|\cdot \|_2$ when $H = I_n$, where $I_n$ is the identity matrix of order $n$. 
The $H$-orthogonality $\bm{x} \perp_H \bm{y}$ represents $(\bm{x}, \bm{y})_H = 0$. 
When $(A\bm{x}, \bm{y})_H = (\bm{x}, A\bm{y})_H$ always satisfies for $\bm{x}, \bm{y} \in \mathbb{R}^{n}$, matrix $A$ is referred to as a self-adjoint in terms of the $H$-inner product. 
This property is equivalent to $(HA)^\top = A^\top H = HA$. 

Following \cite{Sogabe_2005,Kawase_2025}, we also introduce the following for $\hat{\bm{x}}, \hat{\bm{y}} \in \mathbb{R}^{2n}$: 
\begin{align*}
\langle \hat{\bm{x}}, \hat{\bm{y}}\rangle_{\hat H} := \hat{\bm{x}}^\top \hat H\hat{\bm{y}},\quad \hat H := 
\begin{bmatrix}
O & I_n \\
I_n & O
\end{bmatrix} \in \mathbb{R}^{2n\times 2n}. 
\end{align*}
Because $\hat H$ is symmetric but not positive definite, $\langle \cdot, \cdot \rangle_{\hat H}$ is referred to as $\hat H$-quasi-inner product. 
Although the corresponding induced norm cannot be defined, it holds that $\langle \hat{\bm{x}}, \hat{\bm{y}} \rangle_{\hat H} = \langle \hat{\bm{y}}, \hat{\bm{x}} \rangle_{\hat H}$ and $\langle \alpha \hat{\bm{x}} + \hat{\bm{y}}, \hat{\bm{z}} \rangle_{\hat H} = \alpha \langle \hat{\bm{x}}, \hat{\bm{z}} \rangle_{\hat H} + \langle \hat{\bm{y}}, \hat{\bm{z}} \rangle_{\hat H}$ for $\hat{\bm{x}}, \hat{\bm{y}}, \hat{\bm{z}} \in \mathbb{R}^{2n}$ and $\alpha \in \mathbb{R}$. 
For convenience, the $\hat H$-quasi-orthogonality $\hat{\bm{x}} \perp_{\hat H} \hat{\bm{y}}$ represents $\langle \hat{\bm{x}}, \hat{\bm{y}} \rangle_{\hat H} = 0$.

\subsection{Organization}

The remainder of this paper is organized as follows. 
Sect.~\ref{sec2} outlines standard residual smoothing techniques for obtaining the smooth convergence behavior of iterative methods. 
Sect.~\ref{sec3} describes the generation of the Bi-CR residuals using a smoothing scheme. 
After introducing the standard Bi-CG and Bi-CR methods, we present an algorithm for transforming Bi-CG residuals into Bi-CR residuals. 
Sect.~\ref{sec4} presents several lemmas and the main theorem that indicates the bi-orthogonal properties of the transformation algorithm. 
Sect.~\ref{sec5} presents a concise transformation algorithm and numerical example to demonstrate the validity of our insights. 
Finally, Sect.~\ref{sec6} presents the concluding remarks.

\section{Residual smoothing techniques}\label{sec2}

Let $\bm{x}_k$ and $\bm{r}_k (= \bm{b}-A\bm{x}_k)$ be the primary approximations and corresponding residuals, respectively, obtained using an iterative method, where $k$ is the number of iterations. 
Subsequently, residual smoothing techniques generate alternative approximations $\bm{y}_k$ and corresponding residuals $\bm{s}_k (= \bm{b}-A\bm{y}_k)$ using the following forms \cite{Schonauer_1987,Weiss_1996}: 
\begin{align}
&\bm{y}_k = \bm{y}_{k-1} + \eta_k(\bm{x}_k - \bm{y}_{k-1}), \label{smo_x}\\
&\bm{s}_k = \bm{s}_{k-1} + \eta_k(\bm{r}_k - \bm{s}_{k-1}), \label{smo_r}
\end{align}
where $\bm{y}_0 := \bm{x}_0$ and $\bm{s}_0 := \bm{r}_0$, and $\eta_k \in \mathbb{R}$ is a parameter. 
The smooth convergence behavior of residual norms $\|\bm{s}_k\|_H$ can be obtained by appropriately determining $\eta_k$. 
The most typical choice of $\eta_k$ is to satisfy the $H$-orthogonality $\bm{s}_k \perp_H (\bm{r}_k - \bm{s}_{k-1})$ (or to locally minimize $\|\bm{s}_k\|_H$), and $\eta_k$ is given as follows: 
\begin{align}
\eta_k = -\frac{(\bm{s}_{k-1}, \bm{r}_{k} - \bm{s}_{k-1})_H}{(\bm{r}_{k} - \bm{s}_{k-1}, \bm{r}_{k} - \bm{s}_{k-1})_H}. \label{MRS_eta} 
\end{align}
The smoothed residuals then satisfy inequality $\|\bm{s}_k\|_H \leq \min (\|\bm{r}_k\|_H, \|\bm{s}_{k-1}\|_H)$, resulting in a smooth convergence behavior. 
This technique is referred to as MRS. 
The reverse scheme of MRS \cite{Gutknecht_1993} (called orthogonal smoothing \cite{Weiss_1996}), which generates $\bm{r}_k$ using the sequence of $\bm{s}_k$, can also be considered. 

The recursion formula \eqref{smo_r} (and \eqref{smo_x}) can also be viewed as a transformation of the residuals (and approximations) between the different iterative methods. 
For example, CG residual $\bm{r}_k^{\rm cg}$ is transformed into CR residual $\bm{r}_k^{\rm cr}$ using the form \eqref{smo_r} with \eqref{MRS_eta} \cite{Walker_1995,Weiss_1996}. 
Specifically, it holds that 
\begin{align}
\bm{r}_{k}^{\rm cr} = \bm{r}_{k-1}^{\rm cr} + \eta_k(\bm{r}_k^{\rm cg} - \bm{r}_{k-1}^{\rm cr}),\quad \eta_k = -\frac{(\bm{r}_{k-1}^{\rm cr}, \bm{r}_{k}^{\rm cg} - \bm{r}_{k-1}^{\rm cr})}{(\bm{r}_{k}^{\rm cg} - \bm{r}_{k-1}^{\rm cr}, \bm{r}_{k}^{\rm cg} - \bm{r}_{k-1}^{\rm cr})}, \label{cg_cr}
\end{align} 
where $H = I_n$. 
Thus, MRS transforms the residual orthogonal-type methods (such as CG) into residual minimization-type methods (such as CR). 
Furthermore, Bi-CG residual $\bm{r}_k^{\rm bicg}$ is transformed into QMR residual $\bm{r}_k^{\rm qmr}$ using an alternative parameter as follows \cite{Zhou_1994}: 
\begin{align}
\bm{r}_{k}^{\rm qmr} = \bm{r}_{k-1}^{\rm qmr} + \eta_k(\bm{r}_k^{\rm bicg} - \bm{r}_{k-1}^{\rm qmr}),\quad \eta_k = \frac{\tau_k^2}{\rho_k^2}, \label{bicg_qmr}
\end{align}
where $\rho_k^2 := (\bm{r}_k^{\rm bicg}, \bm{r}_k^{\rm bicg})$ and $\tau_k^{-2} = \tau_{k-1}^{-2} + \rho_k^{-2}$ with $\tau_0^2 := \rho_0^2$. 
Note that QMRS is a generalization of the relationship \eqref{bicg_qmr} that can be applied to any sequence of residuals. 

Because the Bi-CG and Bi-CR methods are extensions of the CG and CR methods, respectively, it is natural to wonder whether a transformation exists between the Bi-CG and Bi-CR residuals via the smoothing form~\eqref{smo_r}. 
That is, this study considers the construction of $\eta_k$ such that the following residual transformation holds: 
\begin{align}
\bm{r}_k^{\rm bicr} = \bm{r}_{k-1}^{\rm bicr} + \eta_k (\bm{r}_k^{\rm bicg} - \bm{r}_{k-1}^{\rm bicr}), \label{bicg_bicr}
\end{align}
where $\bm{r}_k^{\rm bicr}$ denotes the Bi-CR residual. 
Note that the Bi-CG method combined with MRS is known as the BICO method and differs from \eqref{bicg_qmr} and \eqref{bicg_bicr}. 
In the next section, we describe the construction of transformation \eqref{bicg_bicr} based on the discussions in \cite{Kawase_2025}.

\section{Generation of the Bi-CR residuals}\label{sec3}

Here, we present the algorithms for the Bi-CG and Bi-CR methods, as well as their residual transformation algorithm. 

\begin{algorithm}[t]
\caption{Standard Bi-CG method \cite{Fletcher_1976}.}\label{Bi-CG}
\begin{algorithmic}[1]
\STATE Select an initial guess $\bm{x}_{0}$.
\STATE Compute $\bm{r}_{0} = \bm{b} - A\bm{x}_{0}$, and choose $\tilde{\bm{r}}_{0}$.
\STATE Set $\bm{p}_{0} = \bm{r}_{0}$ and $\tilde{\bm{p}}_{0} = \tilde{\bm{r}}_{0}$. 
\FOR {$k = 0,1,\dots$, until convergence}
	\STATE $\alpha_{k} = \dfrac{(\tilde{\bm{r}}_{k}, \bm{r}_{k})}{(\tilde{\bm{p}}_{k}, A\bm{p}_{k})}$
	\STATE $\bm{x}_{k+1} = \bm{x}_{k} + \alpha_{k}\bm{p}_{k}$
	\STATE $\bm{r}_{k+1} = \bm{r}_{k} - \alpha_{k}A\bm{p}_{k},\quad \tilde{\bm{r}}_{k+1} = \tilde{\bm{r}}_{k} - \alpha_{k}A^\top \tilde{\bm{p}}_k$
	\STATE $\beta_{k} = \dfrac{(\tilde{\bm{r}}_{k+1}, \bm{r}_{k+1})}{(\tilde{\bm{r}}_{k}, \bm{r}_{k})}$
	\STATE $\bm{p}_{k+1} = \bm{r}_{k+1} + \beta_{k}\bm{p}_{k},\quad \tilde{\bm{p}}_{k+1} = \tilde{\bm{r}}_{k+1} + \beta_{k}\tilde{\bm{p}}_{k}$
\ENDFOR
\end{algorithmic}
\end{algorithm}

\begin{algorithm}[t]
\caption{Standard Bi-CR method \cite{Sogabe_2005,Sogabe_2009}.}\label{Bi-CR}
\begin{algorithmic}[1]
\STATE Select an initial guess $\bm{x}_{0}$.
\STATE Compute $\bm{r}_{0} = \bm{b} - A\bm{x}_{0}$, and choose $\tilde{\bm{r}}_{0}$.
\STATE Set $\bm{p}_{0} = \bm{r}_{0}$ and $\tilde{\bm{p}}_{0} = \tilde{\bm{r}}_{0}$, and compute $\bm{q}_{0} = A\bm{p}_{0}$.
\FOR {$k = 0,1,\dots$, until convergence}
	\STATE $\alpha_{k} = \dfrac{(\tilde{\bm{r}}_{k}, A\bm{r}_{k})}{(A^\top \tilde{\bm{p}}_{k}, \bm{q}_{k})}$
	\STATE $\bm{x}_{k+1} = \bm{x}_{k} + \alpha_{k}\bm{p}_{k}$
	\STATE $\bm{r}_{k+1} = \bm{r}_{k} - \alpha_{k}\bm{q}_{k},\quad \tilde{\bm{r}}_{k+1} = \tilde{\bm{r}}_{k} - \alpha_{k}A^\top \tilde{\bm{p}}_k$
	\STATE $\beta_{k} = \dfrac{(\tilde{\bm{r}}_{k+1}, A\bm{r}_{k+1})}{(\tilde{\bm{r}}_{k}, A\bm{r}_{k})}$
	\STATE $\bm{p}_{k+1} = \bm{r}_{k+1} + \beta_{k}\bm{p}_{k},\quad \tilde{\bm{p}}_{k+1} = \tilde{\bm{r}}_{k+1} + \beta_{k}\tilde{\bm{p}}_{k}$
	\STATE $\bm{q}_{k+1} = A\bm{r}_{k+1} + \beta_{k}\bm{q}_{k}$
\ENDFOR
\end{algorithmic}
\end{algorithm}

\subsection{Bi-CG and Bi-CR algorithms}\label{sec3.1}

The Bi-CG method is typically derived from the Lanczos bi-orthogonalization procedure (e.g., see \cite{Saad_2003}). 
However, the resulting algorithm can also be obtained by applying the CG method to an extended $2n$-dimensional linear system
\begin{align}
\hat A \hat{\bm{x}} = \hat{\bm{b}},\quad \hat A := 
\begin{bmatrix}
A & O \\
O & A^\top
\end{bmatrix},\quad 
\hat{\bm{x}} := 
\begin{bmatrix}
\bm{x} \\
\tilde{\bm{x}}
\end{bmatrix},\quad 
\hat{\bm{b}} := 
\begin{bmatrix}
\bm{b} \\
\tilde{\bm{b}}
\end{bmatrix}. \label{z1}
\end{align}
More precisely, because $\hat H\hat A$ is symmetric, regarding $\hat A$ as self-adjoint in the $\hat H$-quasi-inner product $\langle \cdot, \cdot \rangle_{\hat H}$, one can formally apply the CG algorithm defined in $\langle \cdot, \cdot \rangle_{\hat H}$ to~\eqref{z1}. 
Subsequently, dividing the $2n$-dimensional approximation, residual, and search direction vectors 
\begin{align}
\hat{\bm{x}}_k = 
\begin{bmatrix}
\bm{x}_k \\
\tilde{\bm{x}}_k
\end{bmatrix},\quad 
\hat{\bm{r}}_k = 
\begin{bmatrix}
\bm{r}_k \\
\tilde{\bm{r}}_k
\end{bmatrix},\quad 
\hat{\bm{p}}_k = 
\begin{bmatrix}
\bm{p}_k \\
\tilde{\bm{p}}_k
\end{bmatrix} \label{div_vec}
\end{align}
into their components ($n$-dimensional vectors) yields the standard Bi-CG algorithm shown in Algorithm~\ref{Bi-CG}. 
Here, $\tilde{\bm{r}}_0 \in \mathbb{R}^n$ corresponding to $\tilde{\bm{b}} - A^\top \tilde{\bm{x}}_0$ is chosen such that $(\tilde{\bm{r}}_0, \bm{r}_0) \neq 0$ holds and is referred to as the initial shadow residual. 

The original Bi-CR method is derived in a similar manner. 
That is, applying the CR algorithm defined in $\langle \cdot, \cdot \rangle_{\hat H}$ to \eqref{z1} and reshaping appropriately yields the standard Bi-CR algorithm displayed in Algorithm~\ref{Bi-CR}. 
For details on the derivations above, see~\cite{vanderVorst_2003,Sogabe_2005,Sogabe_2009}.

\subsection{Transformation from Bi-CG residuals into Bi-CR residuals}\label{sec3.2}

As described in \cite{Kawase_2025}, from the derivations of the Bi-CG and Bi-CR algorithms and the transformation \eqref{cg_cr} between the CG and CR residuals, it is naturally considered that the CG method with MRS in $\langle \cdot, \cdot \rangle_{\hat H}$ for \eqref{z1} yields the transformation between the Bi-CG and Bi-CR residuals. 
We briefly review the alternative generation of the Bi-CR residuals using a residual smoothing form. 

Let $\hat{\bm{x}}_k \in \mathbb{R}^{2n}$ and $\hat{\bm{r}}_k (= \hat{\bm{b}} - \hat A\hat{\bm{x}}_k) \in \mathbb{R}^{2n}$ be the primary approximations and the corresponding residuals, respectively, obtained by the CG method in $\langle \cdot, \cdot \rangle_{\hat H}$ for \eqref{z1}. 
Subsequently, based on MRS \eqref{smo_x}--\eqref{MRS_eta}, we generate new approximations $\hat{\bm{y}}_k$ and the corresponding smoothed residuals $\hat{\bm{s}}_k$ in the forms 
\begin{align}
\hat{\bm{y}}_k = \hat{\bm{y}}_{k-1} + \eta_k(\hat{\bm{x}}_k - \hat{\bm{y}}_{k-1}),\quad \hat{\bm{s}}_k = \hat{\bm{s}}_{k-1} + \eta_k(\hat{\bm{r}}_k - \hat{\bm{s}}_{k-1}), \label{hat_y_r}
\end{align}
respectively, where $\hat{\bm{y}}_0 := \hat{\bm{x}}_0$ and $\hat{\bm{s}}_0 := \hat{\bm{r}}_0$, and the parameter $\eta_k \in \mathbb{R}$ is determined to satisfy $\hat H$-quasi-orthogonality $\hat{\bm{s}}_k \perp_{\hat H} (\hat{\bm{r}}_k - \hat{\bm{s}}_{k-1})$ as follows: 
\begin{align}
\eta_k = -\frac{\langle \hat{\bm{s}}_{k-1}, \hat{\bm{r}}_{k} - \hat{\bm{s}}_{k-1}\rangle_{\hat H}}{\langle \hat{\bm{r}}_{k} - \hat{\bm{s}}_{k-1}, \hat{\bm{r}}_{k} - \hat{\bm{s}}_{k-1}\rangle_{\hat H}}. \label{eta_MRS2}
\end{align}
Next, substituting the $2n$-dimensional vectors 
\begin{align*}
\hat{\bm{y}}_k = 
\begin{bmatrix}
\bm{y}_k \\
\tilde{\bm{y}}_k
\end{bmatrix},\quad 
\hat{\bm{s}}_k = 
\begin{bmatrix}
\bm{s}_k \\
\tilde{\bm{s}}_k
\end{bmatrix}
\end{align*}
into \eqref{hat_y_r} and \eqref{eta_MRS2} (together with $\hat{\bm{x}}_k$ and $\hat{\bm{r}}_k$ in \eqref{div_vec}), we obtain the recursion formulas for updating the $n$-dimensional approximations and residuals as follows: 
\begin{align}
\bm{y}_k = \bm{y}_{k-1} + \eta_k(\bm{x}_k - \bm{y}_{k-1}), \quad \bm{s}_k = \bm{s}_{k-1} + \eta_k(\bm{r}_k - \bm{s}_{k-1}), \label{y_s_new} \\
\tilde{\bm{y}}_k = \tilde{\bm{y}}_{k-1} + \eta_k(\tilde{\bm{x}}_k - \tilde{\bm{y}}_{k-1}),\quad \tilde{\bm{s}}_k = \tilde{\bm{s}}_{k-1} + \eta_k(\tilde{\bm{r}}_k - \tilde{\bm{s}}_{k-1}), \label{tilde_y_s_new}
\end{align}
where $\bm{y}_{0} = \bm{x}_{0}$, $\bm{s}_{0} = \bm{r}_{0}$, $\tilde{\bm{y}}_{0} = \tilde{\bm{x}}_{0}$, $\tilde{\bm{s}}_{0} = \tilde{\bm{r}}_{0}$, and 
\begin{align}
\eta_k = -\frac{(\tilde{\bm{s}}_{k-1}, \bm{u}_{k}) + (\bm{s}_{k-1}, \tilde{\bm{u}}_{k})}{2(\tilde{\bm{u}}_{k}, \bm{u}_{k})},\quad \bm{u}_k := \bm{r}_k - \bm{s}_{k-1},\quad \tilde{\bm{u}}_k := \tilde{\bm{r}}_k - \tilde{\bm{s}}_{k-1}. \label{eta_new}
\end{align}
Finally, incorporating \eqref{y_s_new}--\eqref{eta_new} into Algorithm~\ref{Bi-CG} yields Algorithm~\ref{Bi-CR_2}, which represents the transformation from Bi-CG residuals into Bi-CR residuals via a smoothing scheme. 
Here, following \cite{Kawase_2025}, we omit the recursion formula for $\tilde{\bm{y}}_k$. 
Although $\tilde{\bm{s}}_k$ and $\tilde{\bm{u}}_k$ can also be removed by using an alternative form of $\eta_k$ (as described in Sect.~\ref{sec5}), the original algorithm is discussed for convenience. 

The derivation processes indicate that the updating formula for the smoothed residual (i.e., $\bm{s}_k$ in line~11 of Algorithm~\ref{Bi-CR_2} or \eqref{y_s_new}) is equivalent to our target transformation \eqref{bicg_bicr}. 
Numerical experiments in \cite{Kawase_2025} demonstrated that residuals $\bm{s}_k$ generated in Algorithm~\ref{Bi-CR_2} coincide with Bi-CR residual $\bm{r}_k^{\rm bicr}$ generated in Algorithm~\ref{Bi-CR}. 
In the next section, we present the theoretical insights into this equivalence.

\begin{algorithm}[t]
\caption{Residual transformation from Bi-CG into Bi-CR \cite{Kawase_2025}. (Original ver.)}\label{Bi-CR_2}
\begin{algorithmic}[1]
\STATE Select an initial guess $\bm{x}_{0}$.
\STATE Compute $\bm{r}_{0} = \bm{b} -A\bm{x}_{0}$, and choose $\tilde{\bm{r}}_0$.
\STATE Set $\bm{p}_{0} = \bm{r}_{0}$, $\tilde{\bm{p}}_{0} = \tilde{\bm{r}}_{0}$, $\bm{y}_{0} = \bm{x}_{0}$, $\bm{s}_{0} = \bm{r}_{0}$, and $\tilde{\bm{s}}_{0} = \tilde{\bm{r}}_{0}$.
\FOR {$k = 0,1,\dots$, until convergence}
\STATE $\alpha_{k} = \dfrac{(\tilde{\bm{r}}_{k}, \bm{r}_{k})}{(\tilde{\bm{p}}_{k}, A\bm{p}_{k})}$
\STATE $\bm{x}_{k+1} = \bm{x}_{k} + \alpha_{k}\bm{p}_{k}$
\STATE $\bm{r}_{k+1} = \bm{r}_{k} - \alpha_{k}A\bm{p}_{k},\quad \tilde{\bm{r}}_{k+1} = \tilde{\bm{r}}_{k} - \alpha_{k}A^\top \tilde{\bm{p}}_{k}$
\STATE $\bm{u}_{k+1} = \bm{r}_{k+1} - \bm{s}_k,\quad \tilde{\bm{u}}_{k+1} = \tilde{\bm{r}}_{k+1} - \tilde{\bm{s}}_k$
\STATE $\eta_{k+1} = - \dfrac{ (\tilde{\bm{s}}_k, \bm{u}_{k+1} ) + ( \bm{s}_k, \tilde{\bm{u}}_{k+1} ) }{ 2( \tilde{\bm{u}}_{k+1}, \bm{u}_{k+1} ) }$
\STATE $\bm{y}_{k+1} = \bm{y}_k + \eta_{k+1} ( \bm{x}_{k+1} - \bm{y}_k )$
\STATE $\bm{s}_{k+1} = \bm{s}_k + \eta_{k+1} \bm{u}_{k+1},\quad \tilde{\bm{s}}_{k+1} = \tilde{\bm{s}}_k + \eta_{k+1} \tilde{\bm{u}}_{k+1}$
\STATE $\beta_{k} = \dfrac{(\tilde{\bm{r}}_{k+1}, \bm{r}_{k+1})}{(\tilde{\bm{r}}_{k}, \bm{r}_{k})}$
\STATE $\bm{p}_{k+1} = \bm{r}_{k+1} + \beta_{k}\bm{p}_{k},\quad \tilde{\bm{p}}_{k+1} = \tilde{\bm{r}}_{k+1} + \beta_{k}\tilde{\bm{p}}_{k}$
\ENDFOR
\end{algorithmic}
\end{algorithm}

\section{Prove for bi-orthogonal properties}\label{sec4}

As noted in \cite{Sogabe_2005,Sogabe_2009}, from the recursion formulas for updating the residuals, search directions, and their shadow counterparts of the Bi-CR method (i.e., Algorithm~\ref{Bi-CR}), there exist polynomials $R_k(\lambda)$ and $P_k(\lambda)$ of degree $k$ such that 
\begin{align*}
&\bm{r}_k^{\rm bicr} = R_k(A)\bm{r}_0 \in \mathcal{K}_{k+1}(A,\bm{r}_0),\quad \bm{p}_k^{\rm bicr} = P_k(A)\bm{r}_0 \in \mathcal{K}_{k+1}(A,\bm{r}_0),\\
&\tilde{\bm{r}}_k^{\rm bicr} = R_k(A^\top)\tilde{\bm{r}}_0 \in \mathcal{K}_{k+1}(A^\top,\tilde{\bm{r}}_0),\quad \tilde{\bm{p}}_k^{\rm bicr} = P_k(A^\top)\tilde{\bm{r}}_0 \in \mathcal{K}_{k+1}(A^\top,\tilde{\bm{r}}_0) 
\end{align*}
hold, where $R_k(0) = 1$. 
Thus, these vectors are the bases of the corresponding Krylov subspaces. 
The Bi-CR method can then be characterized by the following bi-orthogonal properties \cite{Sogabe_2005,Sogabe_2009}: 
\begin{align*}
&(\bm{r}_i^{\rm bicr}, A^\top \tilde{\bm{r}}_j^{\rm bicr}) = 0 \quad (i \neq j),\\
&(A\bm{p}_i^{\rm bicr}, A^\top \tilde{\bm{p}}_j^{\rm bicr}) = 0 \quad (i \neq j).
\end{align*}
That is, $\bm{r}_k^{\rm bicr}, A\bm{p}_k^{\rm bicr} \perp A^\top \mathcal K_k(A^\top,\tilde{\bm{r}}_0)$ and $\tilde{\bm{r}}_k^{\rm bicr}, A^\top \tilde{\bm{p}}_k^{\rm bicr} \perp A \mathcal K_k(A,\bm{r}_0)$ hold, and these correspond to the well-known bi-orthogonal properties $\bm{r}_k^{\rm bicg}, A\bm{p}_k^{\rm bicg} \perp \mathcal K_k(A^\top,\tilde{\bm{r}}_0)$ and $\tilde{\bm{r}}_k^{\rm bicg}, A^\top \tilde{\bm{p}}_k^{\rm bicg} \perp \mathcal K_k(A,\bm{r}_0)$ in the Bi-CG method. 

Conversely, the recursion formulas for $\bm{s}_k$ and $\tilde{\bm{s}}_k$ in Algorithm~\ref{Bi-CR_2} are similar to those for $\bm{r}_k^{\rm bicr}$ and $\tilde{\bm{r}}_k^{\rm bicr}$, respectively, and auxiliary vectors $\bm{u}_{k+1}$ and $\tilde{\bm{u}}_{k+1}$ can be considered to play the same role of directions $A\bm{p}_k^{\rm bicr}$ and $A^\top \tilde{\bm{p}}_k^{\rm bicr}$, respectively. 
From the updating processes of Algorithm~\ref{Bi-CR_2}, as $\bm{s}_k, \bm{u}_k \in \mathcal K_{k+1}(A,\bm{r}_0)$ and $\tilde{\bm{s}}_k, \tilde{\bm{u}}_k \in \mathcal K_{k+1}(A^\top,\tilde{\bm{r}}_0)$ hold, they can also consist of bases of the corresponding Krylov subspaces. 
Therefore, to confirm the equivalence between $\bm{s}_k$ ($\tilde{\bm{s}}_k$) in Algorithm~\ref{Bi-CR_2} and $\bm{r}_k^{\rm bicr}$ ($\tilde{\bm{r}}_k^{\rm bicr}$), we show that the bi-orthogonal properties 
\begin{align}
&(\bm{s}_i, A^\top \tilde{\bm{s}}_j) = 0 \quad (i \neq j), \label{biortho1} \\
&(\bm{u}_i, \tilde{\bm{u}}_j) = 0 \quad (i \neq j) \label{biortho2}
\end{align}
hold in Algorithm~\ref{Bi-CR_2}. 
This is the goal of this section and the main contribution of this study. 

\subsection{Preparations and observations}

Here, we present preliminary properties to prove \eqref{biortho1} and \eqref{biortho2}. 
Several conclusions for Algorithm~\ref{Bi-CR_2} are presented. 

\begin{proposition}\label{prop}
In Algorithm~\ref{Bi-CR_2}, the $k$th residual $\bm{r}_k$ and search direction $\bm{p}_k$ satisfy the following bi-orthogonal properties: 
\begin{align}
&(\bm{r}_k,\tilde{\bm{r}}_j) = (\bm{r}_k,\tilde{\bm{p}}_j) = (\bm{r}_k,\tilde{\bm{s}}_j) = 0 \quad (j<k), \label{bicg_ortho1} \\ 
&(A\bm{p}_k,\tilde{\bm{r}}_j) = (A\bm{p}_k,\tilde{\bm{p}}_j) = (A\bm{p}_k,\tilde{\bm{s}}_j) = 0 \quad (j<k). \label{bicg_ortho2}
\end{align}
\begin{proof}
Algorithm~\ref{Bi-CR_2} is the same as Algorithm~\ref{Bi-CG} (i.e., the Bi-CG method), except for lines~8--11, which constitute the additional smoothing steps. 
Therefore, $\bm{r}_k$, $\bm{p}_k$, $\tilde{\bm{r}}_k$, and $\tilde{\bm{p}}_k$ are just equivalent to $\bm{r}_k^{\rm bicg}$, $\bm{p}_k^{\rm bicg}$, $\tilde{\bm{r}}_k^{\rm bicg}$, and $\tilde{\bm{p}}_k^{\rm bicg}$, respectively, and thus $\bm{r}_k \perp \tilde{\bm{r}}_j, \tilde{\bm{p}}_j, \tilde{\bm{s}}_j \in \mathcal K_{j+1}(A^\top,\tilde{\bm{r}}_0)$ and $A\bm{p}_k \perp \tilde{\bm{r}}_j, \tilde{\bm{p}}_j, \tilde{\bm{s}}_j \in \mathcal K_{j+1}(A^\top,\tilde{\bm{r}}_0)$ clearly hold for $j < k$ from the bi-orthogonal properties in the Bi-CG method. 
\end{proof}
\end{proposition}

Note that bi-orthogonal properties \eqref{bicg_ortho1} and \eqref{bicg_ortho2} also hold even when the residuals and search directions are reversed to their shadow counterparts; for example, $(\tilde{\bm{r}}_k,\bm{s}_j) = 0$ and $(A^\top \tilde{\bm{p}}_k,\bm{s}_j) = 0$ hold for $j < k$. 

\begin{lemma}\label{lemma1}
Auxiliary vectors $\bm{u}_k := \bm{r}_k - \bm{s}_{k-1}$ and $\tilde{\bm{u}}_k := \tilde{\bm{r}}_k - \tilde{\bm{s}}_{k-1}$ in Algorithm~\ref{Bi-CR_2} satisfy the following recursion formulas for $k \geq 1$. 
\begin{align}
&\bm{u}_{k+1} = (1 - \eta_k)\bm{u}_k - \alpha_kA\bm{p}_k, \label{u_recur} \\
&\tilde{\bm{u}}_{k+1} = (1 - \eta_k)\tilde{\bm{u}}_k - \alpha_kA^{\top}\tilde{\bm{p}}_k. \label{tilde_u_recur}
\end{align}
\begin{proof}
From the recursion formulas of $\bm{r}_k$ and $\bm{s}_k$ in Algorithm~\ref{Bi-CR_2}, \eqref{u_recur} is obtained as 
\begin{align*}
\bm{u}_{k+1} = \bm{r}_{k+1} - \bm{s}_k = \bm{r}_k - \alpha_k A\bm{p}_k - (\bm{s}_{k-1} + \eta_k\bm{u}_k) = (1 - \eta_k)\bm{u}_{k} - \alpha_kA\bm{p}_k. 
\end{align*}
Similarly, \eqref{tilde_u_recur} can be obtained from the recursion formulas of $\tilde{\bm{r}}_k$ and $\tilde{\bm{s}}_k$. 
\end{proof}
\end{lemma}

The above recursive expressions of $\bm{u}_k$ and $\tilde{\bm{u}}_k$ are useful for later proving the bi-orthogonal properties using an induction argument. 

\begin{lemma}\label{lemma2}
In Algorithm~\ref{Bi-CR_2}, the smoothing parameter $\eta_k$ can be rewritten as follows: 
\begin{align}
\eta_k = \frac{(\bm{s}_{k-1}, \tilde{\bm{s}}_{k-1})}{(\bm{u}_k, \tilde{\bm{u}}_k)}, \label{eta_rew}
\end{align}
and the following bi-orthogonality holds for $k \geq 1$: 
\begin{align}
(\bm{s}_k,\tilde{\bm{u}}_k)=0. \label{sk_uk}
\end{align}
\begin{proof}
We have $\bm{r}_k \perp \tilde{\bm{s}}_{k-1}$ (also $\tilde{\bm{r}}_k \perp \bm{s}_{k-1}$) from Proposition~\ref{prop}. 
Thus, with \eqref{eta_new}, it holds that 
\begin{align*}
\eta_k 
= -\frac{(\tilde{\bm{s}}_{k-1}, \bm{r}_k - \bm{s}_{k-1}) + (\bm{s}_{k-1}, \tilde{\bm{r}}_k - \tilde{\bm{s}}_{k-1})}{ 2( \tilde{\bm{u}}_k, \bm{u}_k ) }
= \frac{(\bm{s}_{k-1}, \tilde{\bm{s}}_{k-1})}{(\bm{u}_k, \tilde{\bm{u}}_k)}. 
\end{align*}
Then, by using this form of $\eta_k$ and $\tilde{\bm{r}}_k \perp \bm{s}_{k-1}$, it holds that 
\begin{align*}
(\bm{s}_k,\tilde{\bm{u}}_k) 
&= (\bm{s}_{k-1} + \eta_k\bm{u}_k, \tilde{\bm{u}}_k)
= (\bm{s}_{k-1}, \tilde{\bm{u}}_k) + \frac{(\bm{s}_{k-1},\tilde{\bm{s}}_{k-1})}{(\bm{u}_k,\tilde{\bm{u}}_k)} (\bm{u}_k,\tilde{\bm{u}}_k) \\
&= (\bm{s}_{k-1}, \tilde{\bm{r}}_k - \tilde{\bm{s}}_{k-1}) + (\bm{s}_{k-1}, \tilde{\bm{s}}_{k-1})
=0. 
\end{align*}
\end{proof}
\end{lemma}

From Lemma~\ref{lemma2}, we can determine the difference between conventional MRS and the presented smoothing. 
By simply applying conventional MRS to residual $\bm{r}_k$, the smoothing parameter $\eta_k$ is determined by the orthogonal condition $(\bm{s}_k, \bm{u}_k) = (\bm{s}_k, \bm{r}_k - \bm{s}_{k-1}) = 0$. 
By contrast, the presented smoothing scheme can be interpreted as imposing the bi-orthogonal condition $(\bm{s}_k, \tilde{\bm{u}}_k) = (\bm{s}_k, \tilde{\bm{r}}_k - \tilde{\bm{s}}_{k-1}) = 0$.

\subsection{Main theorem}\label{sec_theorem}

Using the above preparations, we now consider the proof of \eqref{biortho1} and \eqref{biortho2}. 
From the duality of Algorithm~\ref{Bi-CR_2}, we prove the one-sided bi-orthogonality $(\bm{s}_k, A^\top \tilde{\bm{s}}_j) = 0$ for $j < k$, that is sufficient to obtain \eqref{biortho1} (the same holds for \eqref{biortho2}). 
Moreover, because both $\tilde{\bm{s}}_j$ and $\tilde{\bm{p}}_j$ are bases of $\mathcal K_k(A^\top,\tilde{\bm{r}}_0)$ for $j < k$, we may prove $(\bm{s}_k, A^\top \tilde{\bm{p}}_j) = 0$ instead of $(\bm{s}_k, A^\top \tilde{\bm{s}}_j) = 0$. 
These conversions simplify the proof. 

\begin{theorem}\label{theorem}
The following bi-orthogonal properties are satisfied in Algorithm~\ref{Bi-CR_2}: 
\begin{align}
&(\bm{s}_k, A^{\top}\tilde{\bm{p}}_j) = 0 \quad (j<k), \label{biortho3} \\
&(\bm{u}_{k+1}, \tilde{\bm{u}}_{j+1}) = 0 \quad (j<k). \label{biortho4}
\end{align}
\begin{proof}
The proof is obtained by the induction of both the bi-orthogonal properties. 

We first consider the cases where $k=1$ and $j=0$. 
Noting that $\bm{s}_0 = \bm{r}_0$ and $\tilde{\bm{s}}_0 = \tilde{\bm{r}}_0$, smoothing parameter $\eta_1$ in \eqref{eta_rew} can be rewritten as follows: 
\begin{align*}
\eta_1 = \frac{(\bm{s}_0,\tilde{\bm{s}}_0)}{(\bm{u}_1,\tilde{\bm{u}}_1)} = \frac{(\bm{r}_0,\tilde{\bm{r}}_0)}{(\bm{u}_1,\tilde{\bm{r}}_1-\tilde{\bm{s}}_0) } = \frac{\alpha_0 (\bm{p}_0,A^\top \tilde{\bm{p}}_0)}{(\bm{u}_1,\tilde{\bm{r}}_0-\alpha_0A^\top \tilde{\bm{p}}_0 - \tilde{\bm{s}}_0)} = -\frac{(\bm{p}_0,A^\top \tilde{\bm{p}}_0)}{(\bm{u}_1,A^\top \tilde{\bm{p}}_0)}. 
\end{align*}
With $\eta_1$ and $\bm{s}_0 = \bm{p}_0$, we have 
\begin{align}
(\bm{s}_1, A^{\top}\tilde{\bm{p}}_0) = (\bm{s}_0 + \eta_1\bm{u}_1, A^{\top}\tilde{\bm{p}}_0) = (\bm{s}_0, A^{\top}\tilde{\bm{p}}_0) +\eta_1(\bm{u}_1, A^{\top}\tilde{\bm{p}}_0) = 0. \label{bi3_s_k1_j0}
\end{align}
From \eqref{bicg_ortho1}, \eqref{bi3_s_k1_j0}, and $\tilde{\bm{s}}_0 = \tilde{\bm{r}}_0$, we have
\begin{align}
(\bm{u}_2, \tilde{\bm{u}}_1) 
&= (\bm{r}_2 - \bm{s}_1, \tilde{\bm{r}}_1 - \tilde{\bm{s}}_0) = -(\bm{s}_1, \tilde{\bm{r}}_1) +(\bm{s}_1, \tilde{\bm{s}}_0) \notag \\
&= -(\bm{s}_1, \tilde{\bm{r}}_0 - \alpha_0A^\top \tilde{\bm{p}}_0) + (\bm{s}_1, \tilde{\bm{s}}_0) = 0. \label{bi3_u_k1_j0}
\end{align}

Next, as an induction hypothesis, we assume that \eqref{biortho3} and \eqref{biortho4} are satisfied for a certain iteration $k>0$. 
In the following, we show that $(\bm{s}_{k+1}, A^{\top}\tilde{\bm{p}}_j) = 0$ and $(\bm{u}_{k+2}, \tilde{\bm{u}}_{j+1}) = 0$ hold for $j < k + 1$. 

Here, we consider the case of $j = k$. 
From $(\tilde{\bm{r}}_{k+1}, \bm{s}_k)=0$, \eqref{tilde_u_recur}, and \eqref{sk_uk}, the numerator of $\eta_{k+1}$ (cf.~\eqref{eta_rew}) can be expressed as 
\begin{align*}
(\bm{s}_k, \tilde{\bm{s}}_k) 
= (\bm{s}_k, \tilde{\bm{r}}_{k+1} - \tilde{\bm{u}}_{k+1}) 
= - (\bm{s}_k, (1-\eta_k)\tilde{\bm{u}}_k - \alpha_k A^\top \tilde{\bm{p}}_k) = \alpha_k (\bm{s}_k, A^\top \tilde{\bm{p}}_k). 
\end{align*}
The corresponding denominator is expressed as 
\begin{align*}
(\bm{u}_{k+1}, \tilde{\bm{u}}_{k+1}) 
&= (\bm{u}_{k+1}, (1-\eta_k)\tilde{\bm{u}}_k - \alpha_k A^\top \tilde{\bm{p}}_k) = - \alpha_k (\bm{u}_{k+1}, A^\top \tilde{\bm{p}}_k), 
\end{align*}
where $(\bm{u}_{k+1}, \tilde{\bm{u}}_k) = 0$ from the hypothesis. 
Thus, we obtain the alternative form of $\eta_{k+1}$ as 
\begin{align}
\eta_{k+1} = -\frac{(\bm{s}_k,A^\top \tilde{\bm{p}}_k)}{(\bm{u}_{k+1},A^\top \tilde{\bm{p}}_k)}, \label{eta_alt}
\end{align}
and using \eqref{eta_alt} leads to 
\begin{align}
(\bm{s}_{k+1}, A^\top \tilde{\bm{p}}_k) = (\bm{s}_k + \eta_{k+1}\bm{u}_{k+1}, A^\top \tilde{\bm{p}}_k) = (\bm{s}_k, A^\top \tilde{\bm{p}}_k) + \eta_{k+1} (\bm{u}_{k+1}, A^\top \tilde{\bm{p}}_k) = 0. \label{bi_s_k+1_jk}
\end{align}
From \eqref{sk_uk} and \eqref{bicg_ortho1}, we have that 
\begin{align}
(\bm{u}_{k+2}, \tilde{\bm{u}}_{k+1}) = (\bm{r}_{k+2} - \bm{s}_{k+1}, \tilde{\bm{u}}_{k+1}) = (\bm{r}_{k+2}, \tilde{\bm{r}}_{k+1} - \tilde{\bm{s}}_k) = 0. \label{bi_u_k+1_jk}
\end{align}

For $j < k$, using the hypohesis, \eqref{bicg_ortho1}, and \eqref{bicg_ortho2}, we see that 
\begin{align}
(\bm{s}_{k+1}, A^\top \tilde{\bm{p}}_j) 
&= (\bm{s}_k + \eta_{k+1}\bm{u}_{k+1}, A^\top \tilde{\bm{p}}_j) 
= \eta_{k+1}(\bm{u}_{k+1}, A^\top \tilde{\bm{p}}_j) \notag \\
&= \eta_{k+1}(\bm{r}_{k+1} - \bm{s}_k, A^\top \tilde{\bm{p}}_j) 
= 0 \label{bi_s_k+1_j}
\end{align}
and
\begin{align}
(\bm{u}_{k+2}, \tilde{\bm{u}}_{j+1}) 
&= ((1 - \eta_{k+1})\bm{u}_{k+1} - \alpha_{k+1}A\bm{p}_{k+1}, \tilde{\bm{u}}_{j+1}) \notag \\
&= -\alpha_{k+1}(A\bm{p}_{k+1}, \tilde{\bm{u}}_{j+1}) 
= -\alpha_{k+1}(A\bm{p}_{k+1}, \tilde{\bm{r}}_{j+1} - \tilde{\bm{s}}_j) 
= 0. \label{bi_u_k+1_j}
\end{align}
The proof is obtained using \eqref{bi3_s_k1_j0}, \eqref{bi3_u_k1_j0}, and \eqref{bi_s_k+1_jk}--\eqref{bi_u_k+1_j}. 
\end{proof}
\end{theorem}

Consequently, because $\tilde{\bm{p}}_j$ are bases of $\mathcal K_k(A^\top, \tilde{\bm{r}}_0)$ for $j<k$, \eqref{biortho3} leads to $(\bm{s}_k, A^\top \tilde{\bm{s}}_j) = 0\ (j<k)$, and results in \eqref{biortho1} and \eqref{biortho2} from the duality of Algorithm~\ref{Bi-CR_2}.

\section{Concise algorithm and numerical example}\label{sec5}

As \eqref{eta_alt} shows, the smoothing parameter $\eta_{k+1}$ can be represented without using $\tilde{\bm{s}}_k$ and $\tilde{\bm{u}}_{k+1}$. 
When this form of $\eta_{k+1}$, the updates of $\tilde{\bm{s}}_k$ and $\tilde{\bm{u}}_k$ are no longer necessary in Algorithm~\ref{Bi-CR_2}. 
This implies that the Bi-CR residuals can be generated by applying the smoothing forms \eqref{smo_x} and \eqref{smo_r} to the Bi-CG method (i.e., Algorithm~\ref{Bi-CG}) without additional vector updates. 
Algorithm~\ref{Bi-CR_3} represents the concise transformation algorithm, in which $\bm{u}_{k+1}$ is replaced by $\bm{r}_{k+1} - \bm{s}_k$ for readability. 
Note that, in Algorithm~\ref{Bi-CR_3}, $\bm{r}_k$ is Bi-CG residual $\bm{r}_k^{\rm bicg}$ itself and $\bm{s}_k$ coincides with Bi-CR residual $\bm{r}_k^{\rm bicr}$. 

\begin{algorithm}[t]
\caption{Residual transformation from Bi-CG into Bi-CR. (Concise ver.)}\label{Bi-CR_3}
\begin{algorithmic}[1]
\STATE Select an initial guess $\bm{x}_{0}$.
\STATE Compute $\bm{r}_{0} = \bm{b} -A\bm{x}_{0}$, and choose $\tilde{\bm{r}}_0$.
\STATE Set $\bm{p}_{0} = \bm{r}_{0}$, $\tilde{\bm{p}}_{0} = \tilde{\bm{r}}_{0}$, $\bm{y}_{0} = \bm{x}_{0}$, and $\bm{s}_{0} = \bm{r}_{0}$.
\FOR {$k = 0,1,\dots$, until convergence}
\STATE $\alpha_{k} = \dfrac{(\tilde{\bm{r}}_{k}, \bm{r}_{k})}{(\tilde{\bm{p}}_{k}, A\bm{p}_{k})}$
\STATE $\bm{x}_{k+1} = \bm{x}_{k} + \alpha_{k}\bm{p}_{k}$
\STATE $\bm{r}_{k+1} = \bm{r}_{k} - \alpha_{k}A\bm{p}_{k},\quad \tilde{\bm{r}}_{k+1} = \tilde{\bm{r}}_{k} - \alpha_{k}A^\top \tilde{\bm{p}}_{k}$
\STATE $\eta_{k+1} = -\dfrac{(\bm{s}_k,A^\top \tilde{\bm{p}}_k)}{(\bm{r}_{k+1} - \bm{s}_k, A^\top \tilde{\bm{p}}_k)}$
\STATE $\bm{y}_{k+1} = \bm{y}_k + \eta_{k+1} ( \bm{x}_{k+1} - \bm{y}_k )$
\STATE $\bm{s}_{k+1} = \bm{s}_k + \eta_{k+1} ( \bm{r}_{k+1} - \bm{s}_k )$
\STATE $\beta_{k} = \dfrac{(\tilde{\bm{r}}_{k+1}, \bm{r}_{k+1})}{(\tilde{\bm{r}}_{k}, \bm{r}_{k})}$
\STATE $\bm{p}_{k+1} = \bm{r}_{k+1} + \beta_{k}\bm{p}_{k},\quad \tilde{\bm{p}}_{k+1} = \tilde{\bm{r}}_{k+1} + \beta_{k}\tilde{\bm{p}}_{k}$
\ENDFOR
\end{algorithmic}
\end{algorithm}

\begin{figure}[t]
\centering
\includegraphics[scale=0.26]{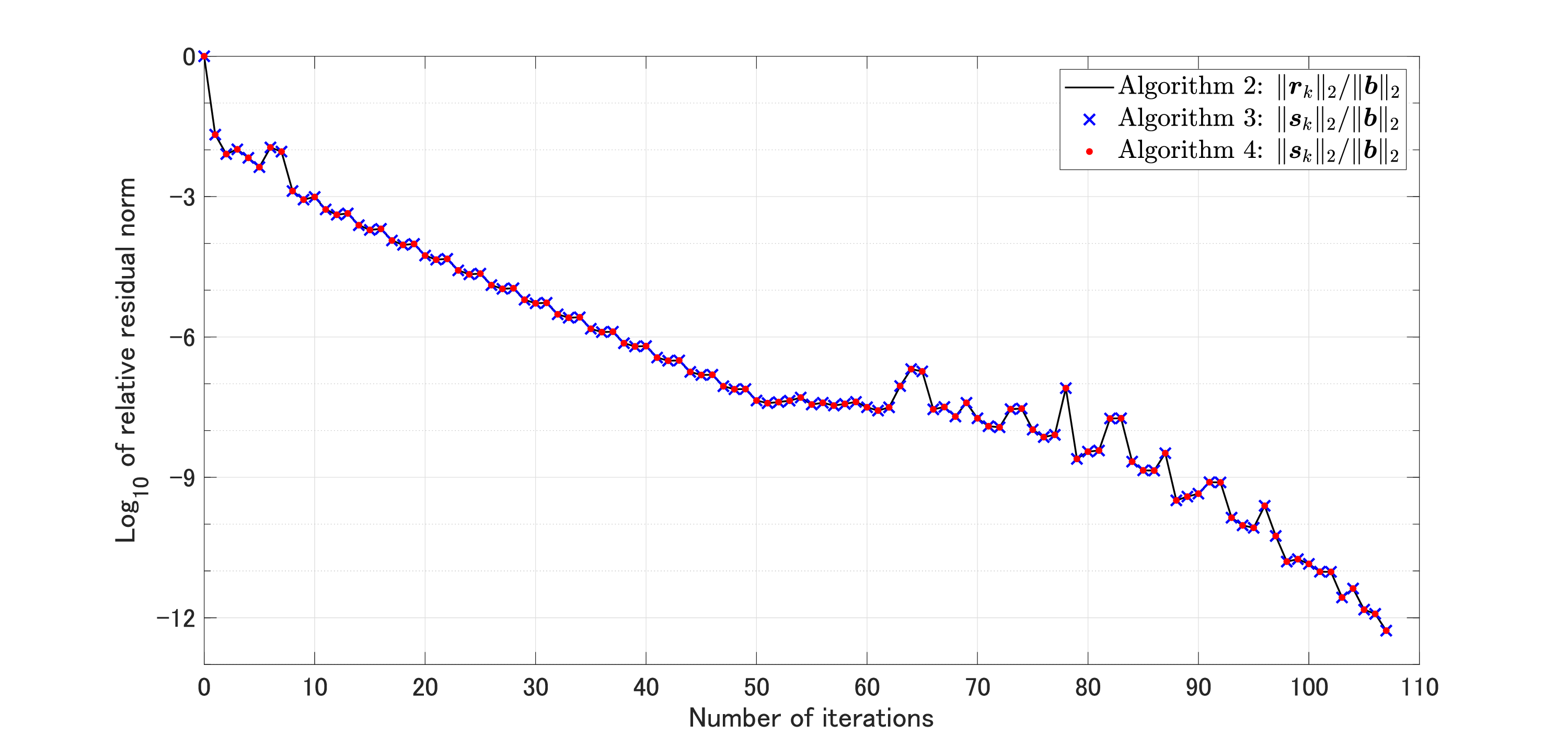}
\caption{Convergence histories of the relative residual norms of Algorithms~\ref{Bi-CR}--\ref{Bi-CR_3} for the Toeplitz matrix.}\label{fig}
\end{figure}

Figure~\ref{fig} shows an example of the convergence histories of Algorithms~\ref{Bi-CR}--\ref{Bi-CR_3} implemented with MATLAB R2021a on a PC (Intel Core i7-1185G7 CPU, 32 GB RAM, double-precision arithmetic). 
The number of iterations and $\log_{10}$ of the relative residual 2-norm are plotted on the horizontal and vertical axes, respectively. 
Following \cite{Sogabe_2005,Kawase_2025}, the Toeplitz matrix 
\begin{align*}
A := 
\begin{bmatrix}
2 & 1 & & & \\
0 & 2 & 1 &  &  \\
1.2 & 0 & 2 & 1 & \\
 & 1.2 & 0 & 2 & \ddots \\
 &  & \ddots & \ddots & \ddots \\
\end{bmatrix}
\in \mathbb{R}^{200\times 200} 
\end{align*}
was used as a test matrix. 
The exact solution and right-hand side were given by $\bm{x}^* := [1,1,\dots,1]^\top$ and $\bm{b} := A\bm{x}^*$, respectively. 
The initial guess and initial shadow residual were set to $\bm{x}_0 := \bm{0}$ and $\tilde{\bm{r}}_0 := \bm{b}$, respectively. 
The iterations were stopped when the relative residual 2-norms were less than $10^{-12}$. 

As shown in Figure~\ref{fig}, the convergence behaviors of the residual norms generated in Algorithms~\ref{Bi-CR}--\ref{Bi-CR_3} coincide, thereby confirming that the Bi-CR residuals can be generated by applying the smoothing scheme (i.e., \eqref{smo_x} and \eqref{smo_r}) to the Bi-CG residuals.

{\remark
We remark on the breakdown of the presented algorithms briefly. 
The conditions for division by zero at $\alpha_k$ and $\beta_k$ in Algorithms~\ref{Bi-CR_2} and \ref{Bi-CR_3} are the same as those in the Bi-CG method (i.e., Algorithm~\ref{Bi-CG}), because these algorithms are identical except for the smoothing steps. 
However, the breakdown conditions for $\eta_k$ in Algorithms~\ref{Bi-CR_2} and \ref{Bi-CR_3} are not entirely clear at present. 
From the discussions in Sect.~\ref{sec4}, the auxiliary vectors $\bm{u}_{k+1}$ and $\tilde{\bm{u}}_{k+1}$ should be collinear to the direction vectors $A\bm{p}_k^{\rm bicr}$ and $A^\top \tilde{\bm{p}}_k^{\rm bicr}$, respectively. 
Therefore, at least if $(A^\top \tilde{\bm{p}}_k^{\rm bicr}, A\bm{p}_k^{\rm bicr})$ (corresponding to the denominator of $\alpha_k$ in Algorihtm~\ref{Bi-CR}) becomes zero, the breakdown also occurs in $\eta_k$ of Algorithm~\ref{Bi-CR_2} (and thus Algorithm~\ref{Bi-CR_3}). 
We leave a further analysis of detailed breakdown conditions open for future work. 
}

\section{Concluding remarks}\label{sec6}

We presented the residual transformation from Bi-CG to Bi-CR using a residual smoothing technique. 
Although the original idea was presented in \cite{Kawase_2025}, the main contribution of this study is a theoretical analysis of the bi-orthogonal properties of the residual and direction vectors in the transformation algorithm (i.e., Algorithm~\ref{Bi-CR_2}). 
Additionally, we derived a more concise algorithm (i.e., Algorithm~\ref{Bi-CR_3}) and generated Bi-CR residuals using this algorithm. 
These observations complement the previous results in \cite{Kawase_2025} and would enable in further advancing the Bi-CG and Bi-CR type methods. 

One potential application of Algorithm~\ref{Bi-CR_3} is to incorporate it with the hybrid procedure~\cite{Brezinski_1994}, which combines different residual sequences for enhancing the convergence. 
This approach may enable us to generate a new residual sequence with superior convergence by combining the Bi-CG and Bi-CR residuals in Algorithm~\ref{Bi-CR_3}, at minimal additional costs. 
A detailed discussion on this regard will be provided in future studies. 
Moreover, we will further analyze the transformation algorithm; for example, the breakdown of the algorithm, numerical behavior in finite precision arithmetic, and similar transformations for the product-type methods of Bi-CG and Bi-CR.

\section*{Acknowledgements}
This study was partly supported by Grant number JP24K14985 from the Grants-in-Aid for Scientific Research Program (KAKENHI) of the Japan Society for the Promotion of Science (JSPS). 
The first author conducted this study during her master's program. 

\section*{Declarations}
The authors declare no conflicts of interest regarding this work.

\end{document}